\documentclass{amsart}

\newtheorem{thm}{Theorem}[section]
\newtheorem{cor}[thm]{Corollary}

\newtheorem{prop}[thm]{Proposition}
\theoremstyle{definition}
\newtheorem{defn}[thm]{Definition}
\theoremstyle{remark}
\newtheorem{rem}[thm]{Remark}

\begin{document}

\title[Foliations on quaternion
CR-submanifolds] {Foliations on quaternion CR-submanifolds}

\author[S. Ianu\c s, A.M Ionescu, G. E. V\^\i lcu]{Stere Ianu\c s,
Adrian Mihai Ionescu, Gabriel Eduard V\^\i lcu}

%\thanks{This work was financially supported by }

\begin{abstract}
The purpose of this paper is to study the canonical foliations of a
quaternion CR-submanifold of a quaternion K\"{a}hler manifold.\\
{\em AMS Mathematics Subject Classification:} 53C15. \\
{\em Key Words and Phrases:} quaternion CR-submanifold, quaternion
K\"{a}hler manifold, foliation.
\end{abstract}

\maketitle

\section{Introduction}

CR-submanifolds of K\"{a}hler manifolds were introduced in
\cite{BJC} by Bejancu. They appear as generalization both of totally
real and of holomorphic submanifolds of K\"{a}hler manifolds.

This notion was further extended in the quaternion settings by
Barros, Chen and Urbano (\cite{BCU}). They consider CR-quaternion
submanifolds of quaternion K\"{a}hlerian manifolds as
generalizations of both quaternion and totally real submanifolds.

If $M$ is a quaternion CR-submanifold  of  a quaternion K\"{a}hler
$\overline{M}$, then two distributions, denoted by $D$ and
$D^\perp$, are defined on $M$. It follows that $D^\perp$ is always
integrable (and of constant rank) i.e. is tangent to a foliation on
$M$. Necessary and sufficient conditions are provided in order that
this foliation became totally geodesic and Riemannian, respectively.

The paper is organized as follows: in Section 2 one reminds basic
definitions and fundamental properties of quaternion CR-submanifolds
of quaternion K\"{a}hler manifolds.

Section 3 allows techniques that will be proven to be useful to
characterize the geometry of $D$ and $D^\perp$; conditions on total
geodesicity are derived.

Section 4 deals with CR-quaternion submanifolds that are ruled by
respect to $D^\perp.$ Section 5 studies the case where $D^\perp$ is
tangent to a Riemannian foliation.

In the last section, QR-products in quaternion space forms are
studied. A characterization of such submanifolds is given.

\section{Preliminaries on quaternion CR-submanifolds}

Let $\overline{M}$ be a differentiable manifold of dimension $n$ and
assume that there is a rank 3-subbundle $\sigma$ of
$End(T\overline{M})$ such that a local basis
$\lbrace{J_1,J_2,J_3}\rbrace$ exists of sections of $\sigma$
satisfying:
       \begin{equation}\label{1}
       \left\{\begin{array}{rcl}
       J_1^2=J_2^2=J_3^2=-Id,\\
       J_1J_2=-J_2J_1=J_3
       \end{array}\right.
       \end{equation}

Then the bundle $\sigma$ is called an almost quaternion structure on
$\overline{M}$ and $\lbrace{J_1,J_2,J_3}\rbrace$ is called canonical
local basis of $\sigma$. Moreover, $\overline{M}$ is said to be an
almost quaternion manifold. It is easy to see that any almost
quaternion manifold is of dimension $n=4m$.

A Riemannian metric $\overline{g}$ is said to be adapted to the
almost quaternion structure $\sigma$ if it satisfies:
 \begin{equation}\label{2}
         \overline{g}(J_\alpha X,J_\alpha Y)=
         \overline{g}(X,Y),\forall\alpha\in\{1,2,3\},
         \end{equation}
for all vector fields $X$,$Y$ on $\overline{M}$ and any local basis
$\lbrace{J_1,J_2,J_3}\rbrace$ of $\sigma$. Moreover,
$(\overline{M},\sigma,\overline{g})$ is said to be an almost
quaternion hermitian manifold.

If the bundle $\sigma$ is parallel with respect to the Levi-Civita
connection $\overline{\nabla}$ of $\overline{g}$, then
$(\overline{M},\sigma,\overline{g})$ is said to be a quaternion
K\"{a}hler manifold. Equivalently, locally defined 1-forms
$\omega_1,\omega_2,\omega_3$ exist such that:
    \begin{equation}\label{3}
          \left\{\begin{array}{rcl}
                 \overline{\nabla}_XJ_1=\omega_3(X)J_2-\omega_2(X)J_3
\\
                 \overline{\nabla}_XJ_2=-\omega_3(X)J_1+\omega_1(X)J_3
\\
                 \overline{\nabla}_XJ_3=\omega_2(X)J_1-\omega_1(X)J_2
                 \end{array}\right.
    \end{equation}
for any vector field $X$ on $\overline{M}$.

Let $(\overline{M},\sigma,\overline{g})$ be an almost quaternion
hermitian manifold. If $X\in T_xM, x\in M$, then the 4-plane $Q(X)$
spanned by $\lbrace {X,J_1X,J_2X,J_3X}\rbrace$ is called a
quaternion 4-plane. A 2-plane in $T_pM$ spanned by
$\lbrace{X,Y}\rbrace$ is called half-quaternion if $Q(X)=Q(Y)$.

The sectional curvature for a half-quaternion 2-plane is called
quaternion sectional curvature. A quaternion K\"{a}hler manifold is
a quaternion space form if its quaternion sectional curvatures are
equal to a constant, say $c$. It is well-known that a quaternion
K\"{a}hler manifold $(\overline{M},\sigma,\overline{g})$ is a
quaternion space form (denoted $\overline{M}(c)$) if and only if its
curvature tensor is:
     \begin{eqnarray}\label{4}
      \overline{R}(X,Y)Z&=&\frac{c}{4}\lbrace \overline{g}(Z,Y)X-
       \overline{g}(X,Z)Y+\sum\limits_{\alpha=1}^3
        [\overline{g}(Z,J_\alpha Y)J_\alpha X\nonumber\\
       &&-\overline{g}(Z,J_\alpha X)J_\alpha Y+
        2\overline{g}(X,J_\alpha Y)J_\alpha Z]\rbrace
      \end{eqnarray}
for all vector fields $X,Y,Z$ on $\overline{M}$ and any local basis
$\lbrace{J_1,J_2,J_3}\rbrace$ of $\sigma.$

\begin{rem}
For a submanifold $M$ of a quaternion K\"{a}hler manifold
$(\overline{M},\sigma,\overline{g})$, we denote by $g$ the metric
tensor induced on $M$. If $\nabla$ is the covariant differentiation
induced on $M$, the Gauss and Weingarten formulas are given by:
\begin{equation}\label{5}
       \overline{\nabla}_XY=\nabla_XY+B(X,Y), \forall X,Y \in
\Gamma(TM)
       \end{equation}
and
\begin{equation}\label{6}
       \overline{\nabla}_XN=-A_NX+\nabla_{X}^{\perp}N, \forall X\in
\Gamma(TM), \forall N\in \Gamma(TM^\perp)
       \end{equation}
where $B$ is the second fundamental form of $M$, $\nabla^\perp$ is
the connection on the normal bundle and $A_N$ is the shape operator
of $M$ with respect to $N$. The shape operator $A_N$ is related to
$h$ by:
\begin{equation}\label{7}
       g(A_NX,Y)=\overline{g}(B(X,Y),N),
       \end{equation}
for all $ X,Y\in \Gamma(TM)$ and $N\in \Gamma(TM^\perp)$.

If we denote by $\overline{R}$ and $R$ the curvature tensor fields
of $\overline{\nabla}$ and $\nabla$ we have the Gauss equation:
\begin{eqnarray}\label{8}
\overline{g}(\overline{R}(X,Y)Z,U)=g(R(X,Y)Z,U)+\overline{g}(B(X,Z),B(Y,U))
       -\overline{g}(B(Y,Z),B(X,U)),\nonumber\\
       \end{eqnarray}
for all $X,Y,Z,U\in \Gamma(TM)$.
\end{rem}

A submanifold $M$ of a quaternion K\"{a}hler manifold
$(\overline{M},\sigma,\overline{g})$ is called a quaternion
CR-submanifold if there exists two orthogonal complementary
distributions $D$ and $D^\perp$ on $M$ such that: \\
i. $D$ is invariant under quaternion structure, that is:
\begin{equation}\label{9}
J_\alpha(D_x)\subseteq D_x, \forall x\in M,
\forall\alpha=\overline{1,3};
\end{equation}
ii. $D^\perp$ is totally real, that is:
\begin{equation}\label{10}
J_\alpha(D_x^\perp)\subseteq T_xM^\perp,
\forall\alpha=\overline{1,3},\forall x\in M.
\end{equation}

A submanifold $M$ of a quaternion K\"{a}hler manifold
$(\overline{M},\sigma,\overline{g})$ is a quaternion submanifold
(respectively, a totally real submanifold) if dim $D^\perp=0$
(respectively, dim $D=0$).

We remark that condition ii. above implies that
$J_\alpha(D_x^\perp)$ are in direct sum, for any local basis as in
(1).

\begin{defn} (\cite{BCU}) Let $M$ be a quaternion CR-submanifold of
a quaternion K\"{a}hler manifold
$(\overline{M},\sigma,\overline{g})$. Then $M$ is called a
QR-product if $M$ is locally the Riemannian product of a quaternion
submanifold and a totally real submanifold of $\overline{M}$.
\end{defn}

\begin{rem}
If we denote by $h^\perp$ and $h$ the second fundamental forms of
$D^\perp$ and $D$, then
we have the following two equations (see \cite{BJC2}):\\
i. $D$-Gauss equation:
\begin{eqnarray}\label{11}
       g(R(X,Y)QZ,QU)&=&g(R^D(X,Y)QZ,QU)+g(h(X,QZ),h(Y,QU))\nonumber\\
       &&-g(h(Y,QZ),h(X,QU)),
       \end{eqnarray}
for all $X,Y,Z,U\in \Gamma(TM)$, where $Q$ is the projection
morphism
of $TM$ on $D$;\\
ii. $D^\perp$-Gauss equation:
\begin{eqnarray}\label{12}
g(R(X,Y)Q{^\perp}Z,Q{^\perp}U)&=&g(R^{D^\perp}(X,Y)Q{^\perp}Z,Q{^\perp}U)+g(h^\perp(X,Q{^\perp}Z),h^\perp(Y,Q{^\perp}U))\nonumber\\
       &&-g(h^\perp(Y,Q{^\perp}Z),h^\perp(X,Q{^\perp}U)),
       \end{eqnarray}
for all $X,Y,Z,U\in \Gamma(TM)$, where $Q^\perp$ is the projection
morphism of $TM$ on $D^\perp$.
\end{rem}

Let $M$ be a quaternion CR-submanifold of a quaternion K\"{a}hler
manifold $(\overline{M},\sigma,\overline{g})$. The we say that:
\\
i. $M$ is $D$-geodesic if:
$$B(X,Y)=0, \forall X,Y\in\Gamma(D).$$
ii. $M$ is $D^\perp$-geodesic if:
$$B(X,Y)=0, \forall X,Y\in\Gamma(D^\perp).$$
iii. $M$ is mixed geodesic if:
$$B(X,Y)=0, \forall X\in\Gamma(D), Y\in\Gamma(D^\perp).$$

We recall now the following result which we shall need in the
sequel.

\begin{thm} $($\cite{BCU}$)$\label{2.5}
Let $M$ be a CR-submanifold of a quaternion K\"{a}hler manifold
$(\overline{M},\sigma,\overline{g})$. Then:
\\
i. The totally real distribution $D^\perp$ is integrable.\\
ii. The
quaternion distribution $D$ is integrable if and only if $M$ is
$D$-geodesic.
\end{thm}

A distribution $D$ in a Riemannian manifold is called minimal if the
trace of its second fundamental form vanishes.

We will illustrate here some of the techniques in this paper on the
following proposition (see also \cite{IP}, \cite{KIM2}).
\begin{prop}\label{2.8}
If $M$ is a CR-submanifold of a quaternion K\"{a}hler manifold
$(\overline{M},\sigma,\overline{g})$, then the quaternion
distribution $D$ is minimal.
\end{prop}
\begin{proof}
Take $X\in\Gamma(D)$ and $U\in\Gamma(D^\perp)$. Then we have:
\begin{eqnarray}
       g(\nabla _X X,U)&=&\overline g(\overline \nabla _X
X,U)\nonumber\\
       &=&\overline g(J_\alpha  \overline \nabla _X X,J_\alpha
U)\nonumber\\
       &=&\overline g( - (\overline \nabla _X J_\alpha  )X + \overline
\nabla _X J_\alpha X,J_\alpha  U)\nonumber\\
       &=&\overline{g}(\omega_\beta(X)J_\gamma X-\omega_\gamma(X)J_\beta X
       +\overline{\nabla}_X J_\alpha X,J_\alpha U)\nonumber\\
       &=&\overline g(\overline\nabla _X J_\alpha  X,J_\alpha
U)\nonumber\\
       &=&\overline g(A_{J_\alpha  U} J_\alpha  X,X)
       \end{eqnarray}
and
\[
g(\nabla _{J_\alpha  X} J_\alpha  X,U) = - \overline g(A_{J_\alpha
U} X,J_\alpha  X)= - \overline g(X,A_{J_\alpha U} J_\alpha  X).
\]

Now for the quaternion distribution $D$ one takes local orthonormal
frame in the form $\{e_i,J_1e_i,J_2e_i,J_3e_i\}$ and summing up over
$i$ will give the assertion.
\end{proof}

\section{Totally real foliation on a quaternion CR-submanifold}

Let $M$ be a quaternion CR-submanifold of a quaternion K\"{a}hler
manifold $(\overline{M},\sigma,\overline{g})$. Then we have the
orthogonal decomposition:
$$TM=D\oplus D^\perp.$$

We have also the following orthogonal decomposition:
$$TM^\perp=\mu\oplus\mu^\perp,$$
where $\mu$ is the subbundle of the normal bundle $TM^\perp$ which
is the orthogonal complement of:
$$\mu^\perp=J_1D^\perp\oplus J_2D^\perp\oplus J_3D^\perp.$$

Since the totally real distribution $D^\perp$ of a quaternion
CR-submanifold $M$ of a quaternion K\"{a}hler manifold
$(\overline{M},\sigma,\overline{g})$ is always integrable we
conclude that we have a foliation $\mathfrak{F}^\perp$ on $M$ with
structural distribution $D^\perp$ and transversal distribution $D$.
We say that $\mathfrak{F}^\perp$ is the canonical totally real
foliation on $M$.

\begin{thm}\label{3.1}
Let $\mathfrak{F}^\perp$ be the canonical totally real foliation on
a quaternion CR-submanifold $M$ of a quaternion K\"{a}hler manifold
$(\overline{M},\sigma,\overline{g})$. The next assertions are
equivalent:
\\
i. $\mathfrak{F}^\perp$ is totally geodesic;
\\
ii. $B(X,Y)\in\Gamma(\mu)$, $\forall X\in\Gamma(D)$,
$Y\in\Gamma(D^\perp)$;
\\
iii. $A_NX\in\Gamma(D^\perp)$, $\forall X\in\Gamma(D^\perp)$,
$N\in\Gamma(\mu^\perp)$;
\\
iv. $A_NY\in\Gamma(D)$, $\forall Y\in\Gamma(D)$,
$N\in\Gamma(\mu^\perp)$.
\end{thm}
\begin{proof}
For $X,Z\in\Gamma(D^\perp)$ and $Y\in\Gamma(D)$ we have:
\begin{eqnarray}
\overline{g}(J_\alpha(\nabla_XZ),Y)&=&-\overline{g}(\overline{\nabla}_XZ-B(X,Z),J_\alpha Y)\nonumber\\
&=&\overline{g}(-(\overline{\nabla}_XJ_\alpha)Z+\overline{\nabla}_XJ_\alpha
       Z,Y)\nonumber\\
       &=&\overline{g}(\omega_\beta(X)J_\gamma
Z-\omega_\gamma(X)J_\beta Z+\overline{\nabla}_XJ_\alpha
       Z,Y)\nonumber\\
       &=&\overline{g}(-A_{J_\alpha Z}X+\nabla_X^{\perp}J_\alpha
Z,Y)\nonumber\\
       &=&-g(A_{J_\alpha Z}X,Y)\nonumber
       \end{eqnarray}
where $(\alpha,\beta,\gamma)$ is an even permutation of (1,2,3), and
taking into account of (\ref{7}) we obtain:
\begin{equation}\label{13}
\overline{g}(J_\alpha(\nabla_XZ),Y)=-\overline{g}(B(X,Y),J_\alpha
Z).
\end{equation}
i. $\Rightarrow$ ii.

If $\mathfrak{F}^\perp$ is totally geodesic, then
$\nabla_XZ\in\Gamma(D^\perp)$, for $X,Z\in\Gamma(D^\perp)$ and from
(\ref{13}) we derive:
$$\overline{g}(B(X,Y),J_\alpha
Z)=0$$ and the implication is clear.
\\
ii. $\Rightarrow$ i.

If we suppose $B(X,Y)\in\Gamma(\mu)$, $\forall X\in\Gamma(D)$,
$Y\in\Gamma(D^\perp)$, then from (\ref{13}) we derive:
$$\overline{g}(J_\alpha(\nabla_XZ),Y)=0$$
and we conclude $\nabla_XZ\in\Gamma(D^\perp)$. Thus
$\mathfrak{F}^\perp$ is totally geodesic.
\\
ii. $\Leftrightarrow$ iii.

This equivalence is clear from (\ref{7}).
\\
iii. $\Leftrightarrow$ iv.

This equivalence is true because $A_N$ is a self-adjoint operator.

\end{proof}

\begin{cor}\label{3.3}
Let $\mathfrak{F}^\perp$ be the canonical totally real foliation on
a quaternion CR-submanifold $M$ of a quaternion K\"{a}hler manifold
$(\overline{M},\sigma,\overline{g})$. If $M$ is mixed geodesic, then
$\mathfrak{F}^\perp$ is totally geodesic.
\end{cor}
\begin{proof}
The assertion is clear from Theorem \ref{3.1}.
\end{proof}

\begin{cor}\label{3.4}
Let $\mathfrak{F}^\perp$ be the canonical totally real foliation on
a quaternion CR-submanifold $M$ of a quaternion K\"{a}hler manifold
$(\overline{M},\sigma,\overline{g})$ with $\mu=0$. Then $M$ is mixed
geodesic if and only if $\mathfrak{F}^\perp$ is totally geodesic.
\end{cor}
\begin{proof}
The assertion is immediate from Theorem \ref{3.1}.
\end{proof}

\section{Totally real ruled quaternion CR-submanifolds}

A submanifold $M$ of a Riemannian manifold
$(\overline{M},\overline{g})$ is said to be a ruled submanifold if
it admits a foliation whose leaves are totally geodesic immersed in
$(\overline{M},\overline{g})$.

\begin{defn}
A quaternion CR-submanifold of a quaternion K\"{a}hler manifold
which is a ruled submanifold with respect to the foliation
$\mathfrak{F}^\perp$ is called totally real ruled quaternion
CR-submanifold.
\end{defn}

\begin{thm}\label{4.2}
Let $M$ be a quaternion CR-submanifold of a quaternion K\"{a}hler
manifold $(\overline{M},\sigma,\overline{g})$. The next assertions
are equivalent:
\\
i. $M$ is a totally real ruled quaternion CR-submanifold.
\\
ii. $M$ is $D^\perp$-geodesic and:
$$B(X,Y)\in\Gamma(\mu),\ \forall
X\in\Gamma(D),\ Y\in\Gamma(D^\perp).$$ iii. The subbundle
$\mu^\perp$ is $D^\perp$-parallel, i.e:
$$\nabla_X^\perp J_\alpha Z\in\Gamma(\mu^\perp),\ \forall X,Z\in
D^\perp,\ \alpha=\overline{1,3}$$ and the second fundamental form
satisfies:
$$B(X,Y)\in\Gamma(\mu),\ \forall
X\in\Gamma(D^\perp),\ Y\in\Gamma(TM).$$ iv. The shape operator
satisfies:
$$A_{J_\alpha Z}X=0,\ \forall X,Z\in D^\perp,\ \alpha=\overline{1,3}$$
and
$$A_N X\in\Gamma(D),\ \forall X\in\Gamma(D^\perp),\ N\in\Gamma(\mu).$$
\end{thm}
\begin{proof} i. $\Leftrightarrow$ ii. For any $X,Z\in\Gamma(D^\perp)$ we have:
\begin{eqnarray}
       \overline{\nabla}_XZ&=&\nabla_XZ+B(X,Z)\nonumber\\
       &=&\nabla_X^{D^\perp}Z+h^\perp(X,Z)+B(X,Z)\nonumber
       \end{eqnarray}
and thus we conclude that the leafs of $D^\perp$ are totally
geodesic immersed in $\overline{M}$ iff $h^\perp=0$ and $M$ is
$D^\perp$-geodesic.

The equivalence is now clear from Theorem \ref{3.1}.
\\
i. $\Leftrightarrow$ iii. For any $X,Z\in\Gamma(D^\perp)$, and
$U\in\Gamma(D)$ we have:
\begin{eqnarray}
\overline{g}(\overline{\nabla}_XZ,U)&=&\overline{g}(J_\alpha\overline{\nabla}_XZ,J_\alpha U)\nonumber\\
&=&\overline{g}(-(\overline{\nabla}_XJ_\alpha)Z+\overline{\nabla}_XJ_\alpha
Z,J_\alpha
       U)\nonumber\\
       &=&\overline{g}(\omega_\beta(X)J_\gamma
Z-\omega_\gamma(X)J_\beta Z+\overline{\nabla}_XJ_\alpha
       Z,J_\alpha U)\nonumber\\
       &=&\overline{g}(-A_{J_\alpha Z}X+\nabla_X^{\perp}J_\alpha
Z,J_\alpha U)\nonumber\\
       &=&-g(A_{J_\alpha Z}X,J_\alpha U)\nonumber
       \end{eqnarray}
where $(\alpha,\beta,\gamma)$ is an even permutation of (1,2,3), and
taking into account of (\ref{7}) we obtain:
\begin{equation}\label{14}
\overline{g}(\overline{\nabla}_XZ,U)=-\overline{g}(B(X,J_\alpha
U),J_\alpha Z).
\end{equation}

On the other hand, for any $X,Z,W\in\Gamma(D^\perp)$ we have:
\begin{eqnarray}\label{15}
       \overline{g}(\overline{\nabla}_XZ,J_\alpha
W)&=&\overline{g}(\overline{\nabla}_XZ+B(X,Z),J_\alpha W)\nonumber\\
       &=&\overline{g}(B(X,Z),J_\alpha W).
       \end{eqnarray}

If $X,Z\in\Gamma(D^\perp)$ and $N\in\Gamma(\mu)$, then we have:
\begin{eqnarray}
\overline{g}(\overline{\nabla}_XZ,N)&=&\overline{g}(J_\alpha\overline{\nabla}_XZ,J_\alpha N)\nonumber\\
&=&\overline{g}(-(\overline{\nabla}_XJ_\alpha)Z+\overline{\nabla}_XJ_\alpha Z,J_\alpha N)\nonumber\\
       &=&\overline{g}(\omega_\beta(X)J_\beta
Z+\omega_\gamma(X)J_\gamma
       Z-J_\alpha \overline{\nabla}_XJ_\alpha Z,N)\nonumber\\
       &=&\overline{g}(\overline{\nabla}_XJ_\alpha Z,J_\alpha
N)\nonumber\\
       &=&\overline{g}(-A_{J_\alpha Z}X+\nabla_X^{\perp}J_\alpha
Z,J_\alpha N)\nonumber
       \end{eqnarray}
and thus we obtain:
\begin{equation}\label{16}
\overline{g}(\overline{\nabla}_XZ,N)=\overline{g}(\nabla^\perp_X
J_\alpha Z,J_\alpha N).
\end{equation}

Finally, $M$ is a totally real ruled quaternion CR-submanifold iff
$\overline{\nabla}_X Z\in\Gamma(D^\perp)$, $\forall
X,Z\in\Gamma(D^\perp)$ and by using (\ref{14}), (\ref{15}) and
(\ref{16}) we deduce the equivalence.
\\
ii. $\Leftrightarrow$ iv.

This is clear from (\ref{7}).
\end{proof}

\begin{cor}\label{4.3}
Let $M$ be a quaternion CR-submanifold of a quaternion K\"{a}hler
manifold $(\overline{M},\sigma,\overline{g})$.  If $M$ is totally
geodesic, then $M$ is a totally real ruled quaternion
CR-submanifold.
\end{cor}
\begin{proof}
The assertion is clear from Theorem \ref{4.2}.
\end{proof}

\section{Riemannian foliations and quaternion CR-submanifolds}

Let $(M,g)$ be a Riemannian manifold and $\mathfrak{F}$ a foliation
on $M$. The metric $g$ is said to be bundle-like for the foliation
$\mathfrak{F}$ if the induced metric on $D^\perp$ is parallel with
respect to the intrinsic connection on the transversal distribution
$D^\perp$. This is true if and only if the Levi-Civita connection
$\nabla$ of $(M,g)$ satisfies (see \cite{BJC2}):
\begin{equation}\label{17}
g(\nabla_{Q^\perp Y}QX,Q^\perp Z)+g(\nabla_{Q^\perp Z}QX,Q^\perp
Y)=0,\ \forall X,Y,Z\in\Gamma(TM).
\end{equation}

If for a given foliation $\mathfrak{F}$ there exists a Riemannian
metric $g$ on $M$ which is bundle-like for $\mathfrak{F}$, then we
say that $\mathfrak{F}$ is a Riemannian foliation on $(M,g)$.

\begin{thm}\label{5.1}
Let $M$ be a quaternion CR-submanifold of a quaternion K\"{a}hler
manifold $(\overline{M},\sigma,\overline{g})$. The next assertions
are equivalent:
\\
i. The induced metric $g$ on $M$ is bundle-like for totally real
foliation $\mathfrak{F}^\perp$.
\\
ii. The second fundamental form $B$ of $M$ satisfies:
$$B(U,J_\alpha
V)+B(V,J_\alpha U)\in\Gamma(\mu)\oplus J_\beta(D^\perp)\oplus
J_\gamma(D^\perp),$$ for any $U,V\in\Gamma(D)$ and $\alpha=1,2$ or
$3$, where $(\alpha,\beta,\gamma)$ is an even permutation of
$(1,2,3)$.
\end{thm}
\begin{proof}
From (\ref{17}) we deduce that $g$ is bundle-like for totally real
foliation $\mathfrak{F}^\perp$ iff:
\begin{equation}\label{18}
g(\nabla_U X,V)+g(\nabla_V X,U)=0,\ \forall X\in\Gamma(D^\perp),\
U,V\in\Gamma(D).
\end{equation}

On the other hand, for any $X\in\Gamma(D^\perp),\ U,V\in\Gamma(D)$
we have:
\begin{eqnarray}
       g(\nabla_U X,V)+g(\nabla_V
X,U)&=&\overline{g}(\overline{\nabla}_U X-B(U,X),V)+
       \overline{g}(\overline{\nabla}_V X-B(V,X),U)\nonumber\\
       &=&\overline{g}(\overline{\nabla}_U X,V)+
       \overline{g}(\overline{\nabla}_V X,U)\nonumber\\
       &=&\overline{g}(-(\overline{\nabla}_U
J_\alpha)X+\overline{\nabla}_U J_\alpha X,J_\alpha V)\nonumber\\
       &&+\overline{g}(-(\overline{\nabla}_V
J_\alpha)X+\overline{\nabla}_V J_\alpha X,J_\alpha U)\nonumber\\
       &=&\overline{g}(\omega_\beta(U)J_\gamma
X-\omega_\gamma(U)J_\beta X
       +\overline{\nabla}_U J_\alpha X,J_\alpha V)\nonumber\\
       &&+\overline{g}(\omega_\beta(V)J_\gamma
X-\omega_\gamma(V)J_\beta X
       +\overline{\nabla}_V J_\alpha X,J_\alpha U)\nonumber\\
       &=&\overline{g}(\overline{\nabla}_U J_\alpha X,J_\alpha V)
       +\overline{g}(\overline{\nabla}_V J_\alpha X,J_\alpha
U)\nonumber\\
       &=&-g(A_{J_\alpha X} U,J_\alpha V)-g(A_{J_\alpha X}V, J_\alpha
U)\nonumber
       \end{eqnarray}
where $(\alpha,\beta,\gamma)$ is an even permutation of (1,2,3), and
taking into account of (\ref{7}) we derive:
\begin{equation}\label{19}
g(\nabla_U X,V)+g(\nabla_V X,U)=-\overline{g}(B(U,J_\alpha
V)+B(V,J_\alpha U),J_\alpha X),
\end{equation}
for any $X\in\Gamma(D^\perp),\ U,V\in\Gamma(D)$.

The proof is now complete from (\ref{18}) and (\ref{19}).
\end{proof}

\begin{cor}\label{5.2}
Let $\mathfrak{F}^\perp$ be the canonical totally real foliation on
a quaternion CR-submanifold $M$ of a quaternion K\"{a}hler manifold
$(\overline{M},\sigma,\overline{g})$ with $\mu=0$. Then the induced
metric $g$ on $M$ is bundle-like for $\mathfrak{F}^\perp$ if and
only if the second fundamental form $B$ of $M$ satisfies:
$$B(U,J_\alpha
V)+B(V,J_\alpha U)\in J_\beta(D^\perp)\oplus J_\gamma(D^\perp),$$
for any $U,V\in\Gamma(D)$ and $\alpha=1,2$ or $3$, where
$(\alpha,\beta,\gamma)$ is an even permutation of $(1,2,3)$.
\end{cor}
\begin{proof}
The assertion is immediate from Theorem \ref{5.1}.
\end{proof}

\begin{cor}
Let $\mathfrak{F}^\perp$ be the canonical totally real foliation on
a quaternion CR-submanifold $M$ of a quaternion K\"{a}hler manifold
$(\overline{M},\sigma,\overline{g})$ with $\mu=0$. Then
$\mathfrak{F}^\perp$ is totally geodesic with bundle-like metric
$g=\overline{g}_{|M}$ if and only if $M$ is mixed geodesic and the
second fundamental form $B$ of $M$ satisfies:
$$B(U,J_\alpha
V)+B(V,J_\alpha U)\in J_\beta(D^\perp)\oplus J_\gamma(D^\perp),$$
for any $U,V\in\Gamma(D)$ and $\alpha=1,2$ or $3$, where
$(\alpha,\beta,\gamma)$ is an even permutation of $(1,2,3)$.
\end{cor}
\begin{proof}
The proof follows from Corollary \ref{3.4} and Corollary \ref{5.2}.
\end{proof}

\section{QR-products in quaternion K\"{a}hler manifolds}

From Theorem \ref{2.5} we deduce that any $D$-geodesic
CR-submanifold of a quaternion K\"{a}hler manifold admits a
$\sigma$-invariant totally geodesic foliation, which we denote by
$\mathfrak{F}$.

\begin{prop}\label{6.2}
If $M$ is a totally geodesic quaternion CR-submanifold of a
quaternion K\"{a}hler manifold $(\overline{M},\sigma,\overline{g})$,
then $M$ is a ruled submanifold with respect to both foliations
$\mathfrak{F}$ and $\mathfrak{F^\perp}$.
\end{prop}
\begin{proof}
The assertion follows from Corollary \ref{4.3} and Theorem
\ref{2.5}.
\end{proof}

\begin{thm}\label{6.3}
Let $M$ be a quaternion CR-submanifold of a quaternion K\"{a}hler
manifold $(\overline{M},\sigma,\overline{g})$. Then $M$ is a
QR-product if and only if the next three conditions are satisfied:
\\
i. $M$ is $D$-geodesic.
\\
ii. $M$ is $D^\perp$-geodesic.
\\
iii. $B(X,Y)\in\Gamma(\mu),\ \forall X\in\Gamma(D^\perp),\
Y\in\Gamma(D)$.
\end{thm}
\begin{proof}
The proof is immediate from Theorems \ref{2.5} and \ref{4.2}.
\end{proof}

\begin{cor}
Any totally geodesic quaternion CR-submanifold of a quaternion
K\"{a}hler manifold is a QR-product.
\end{cor}
\begin{proof}
The assertion is clear.
\end{proof}

\begin{cor}
Let $M$ be a quaternion CR-submanifold of a quaternion K\"{a}hler
manifold $(\overline{M},\sigma,\overline{g})$ with $\mu=0$. Then $M$
is a QR-product if and only if $M$ is totally geodesic.
\end{cor}
\begin{proof}
The assertion is immediate from Theorem \ref{6.3}.
\end{proof}

\begin{thm}
Let $M$ be a quaternion CR-submanifold of a quaternion space form
$\overline{M}(c)$. If $M$ is $D$-geodesic and ruled submanifold with
respect to totally real foliation $\mathfrak{F}^\perp$, then $M$ is
a QR-product. Moreover, locally $M$ is a Riemannian product $L\times
L^\perp$, where $L$ is a quaternion space form, having quaternion
sectional curvature $c$, and $L^\perp$ is a real space form, having
sectional curvature $\frac{c}{4}$.
\end{thm}
\begin{proof}

Because $\mathfrak{F}^\perp$ is a totally geodesic foliation, from
(\ref{12}) we obtain:
\begin{equation}\label{20}
g(R(X,Y)Z,U)=g(R^{D^\perp}(X,Y)Z,U)
\end{equation}
for any $X,Y\in\Gamma(TM)$, $Z,U\in\Gamma(D^\perp)$.

From Gauss equation and (\ref{20}) we obtain:
\begin{eqnarray}\label{21}
\overline{g}(\overline{R}(X,Y)Z,U)&=&g(R^{D^\perp}(X,Y)Z,U)+\overline{g}(B(X,Z),B(Y,U))\nonumber\\
       &&-\overline{g}(B(Y,Z),B(X,U)).
       \end{eqnarray}

Since $M$ is a ruled submanifold with respect to totally real
foliation $\mathfrak{F}^\perp$, then $M$ is $D^\perp$-geodesic (by
Theorem \ref{4.2}) and from (\ref{21}) we derive:
\begin{equation}\label{22}
       \overline{g}(\overline{R}(X,Y)Z,U)=g(R^{D^\perp}(X,Y)Z,U)
       \end{equation}
for any $X,Y,Z,U\in\Gamma(D^\perp)$.

Now, from (\ref{4}) and (\ref{22}) we derive:
\begin{equation}\label{24}
g(R^{D^\perp}(X,Y)Y,X)=\overline{g}(\overline{R}(X,Y)Y,X)=\frac{c}{4},
\end{equation}
for any orthogonal unit vector fields $X,Y\in\Gamma(D^\perp)$. Thus
we conclude that the leaves of the totally real foliation
$\mathfrak{F}^\perp$ are of constant curvature $\frac{c}{4}$.

Since $M$ is $D$-geodesic, from Gauss equation we obtain:
\begin{equation}\label{25}
       \overline{g}(\overline{R}(X',Y')Z',U')=g(R(X',Y')Z',U'),
       \end{equation}
for any $X',Y',Z',U'\in\Gamma(D)$.

On the other hand, if $M$ is $D$-geodesic, then $\mathfrak{F}$ is a
totally geodesic foliation and from (\ref{11}) we obtain:
\begin{equation}\label{26}
       g(R(X',Y')Z',U')=g(R^D(X',Y')Z',U'),
       \end{equation}
for any $X',Y'\in\Gamma(TM)$, $Z',U'\in\Gamma(D)$.

From (\ref{25}) and (\ref{26}) we deduce:
\begin{equation}\label{27}
       \overline{g}(\overline{R}(X',Y')Z',U')=g(R^D(X',Y')Z',U'),
       \end{equation}
for any $X',Y',Z',U'\in\Gamma(D)$.

Now, from (\ref{4}) and (\ref{27}) we derive:
\begin{equation}\label{28}
g(R^{D^\perp}(X',J_\alpha X')J_\alpha
X',X')=\overline{g}(\overline{R}(X',J_\alpha X')J_\alpha X',X')=c,
\end{equation}
for any unit vector field $X'\in\Gamma(D)$. Thus we conclude that
the the leaves of the foliation $\mathfrak{F}$ are of constant
quaternion sectional curvature $c$.

The proof is now complete from (\ref{24}) and (\ref{28}).
\end{proof}

\

\begin{center}

Stere Ianu\c s \\
{\em    University of Bucharest,\\
        Department of Mathematics,\\
        C.P. 10-119, Post. Of. 10, Bucharest 72200, Romania}\\
        e-mail: ianus@gta.math.unibuc.ro
\end{center}
\

\begin{center}

Adrian Mihai Ionescu \\
{\em    University Politehnica of Bucuresti,\\
        Department of Mathematics,\\
        Splaiul Independentei, Nr. 313, Sector 6, Bucure\c sti, Romania}\\
        e-mail: aionescu@math.pub.ro
\end{center}
\

\begin{center}
Gabriel Eduard V\^\i lcu \\
{\em  ''Petroleum-Gas'' University of Ploie\c sti,\\
         Department of Mathematics and Computer Science,\\
         Bulevardul Bucure\c sti, Nr. 39, Ploie\c sti, Romania}\\
         e-mail: gvilcu@mail.upg-ploiesti.ro\\
\end{center}

\end{document}